\documentclass{amsart}

\usepackage{amsmath}
\usepackage{amscd}
\usepackage{amssymb} 

\newcommand{\cal}{\mathcal}

\newcommand{\bZ}{{\Bbb Z}}
\newcommand{\cA}{{\cal A}}

\newcommand{\cH}{{\cal H}}

\newcommand{\fd}{{\frak d}}

\DeclareMathOperator{\End}{End}

\DeclareMathOperator{\Img}{Im}

\DeclareMathOperator{\Ker}{Ker}

\newtheorem{theorem}{Theorem}[section]
\newtheorem{theorem/definition}{Theorem/Definition}[section]

\newtheorem{lemma}{Lemma}[section]

\theoremstyle{remark}

\theoremstyle{definition}

\begin{document}
\title
{Hodge theory and $A_{\infty}$ structures on cohomology}
\author{Jian Zhou}
\address{Department of Mathematics\\
Texas A\&M University\\
College Station, TX 77843}
\email{zhou@math.tamu.edu}
\begin{abstract}
We use Hodge theory and a construction of Merkulov
to construct $A_{\infty}$ structures on
de Rham cohomology and Dolbeault cohomology.  
\end{abstract}
\maketitle

Hodge theory is a powerful tool in differential geometry.
Classically,  it can be used to identify the de Rham cohomology
of a closed oriented Riemannian manifold with the space of 
harmonic forms on it as vector spaces.
The wedge product on differential forms provides an
algebra structure on the de Rham cohomology.
By virtue of being isomorphic to the de Rham cohomology,
the space of harmonic forms has an induced associative multiplication.
However, the wedge product of two harmonic forms may not be harmonic.
One needs to define a multiplication of two harmonic forms 
by taking the harmonic part of their wedge product,
and then show that this multiplication is indeed associative 
and can be identified with the wedge product on the de Rham cohomology.
In this paper, 
we show that one can actually take advantage of this awkward situation
to construct higher multiplications and define
a structure of $A_{\infty}$ algebra on the space of harmonic forms.

Originally,
$A_{\infty}$ structures were introduced by Stasheff \cite{Sta1, Sta2}
in 1963 in the study of $H$ spaces.
Together with various cousins, 
they have appeared in the last couple of decades 
in many places in Mathematics and Mathematical Physics.
In particular,
infinite structures are very useful in the formulation of mirror symmetry.
For example, Fukaya \cite{Fuk} constructed $A_{\infty}$ categories from
symplectic manifolds, 
which is used in the formulation of homological mirror symmetry 
by Kontsevich \cite{Kon}.
For Calabi-Yau manifold,
where the notion of mirror symmetry was 
originally conceived \cite{Yau, Gre-Yau},
recent formulations of the mirror symmetry use the notion of
Frobenius manifolds introduced by Dubrovin \cite{Dub}.
As pointed out in Manin \cite{Man1},
a formal Frobenius manifold structure on a vector space 
with a nondegenerate pairing
is equivalent to a cyclic $Comm_{\infty}$-structure on it.
The popular theory of quantum cohomology provides 
construction of formal Frobenius manifold structures 
on de Rham cohomology of symplectic manifolds. 
The appearance of infinite algebra structures in 
the theory of quantum cohomology and mirror symmetry
indicates the importance of 
the study of infinite algebras in differential geometry.
 
Our construction is based on  a recent paper of Merkulov \cite{Mer},
where he gave  a nice construction of $A_{\infty}$ algebra.
Together with Hodge theory of closed K\"{a}hler manifold,
he constructed an $A_{\infty}$ structure on
a subcomplex of the de Rham complex which contains the harmonic forms.
Similar constructions can be carried out for the deformation complex of
Calabi-Yau manifolds.
The simple observation of this paper is that Merkulov's construction
can be used to give a construction of
$A_{\infty}$ struction on the space of harmonic forms
of any oriented closed Riemannian manifold.
Similarly,
Hodge theory of the deformation complex of any closed complex $M$
manifold can be used to construct an $A_{\infty}$ structure on
$H^{-*, *}(M)$.
Similar constructions can be carried out for Dolbeault cohomology 
of endomorphism bundle of a holomorphic vector bundle over a complex manifold.

The rest of the paper is arranged as follows.
We review the definition of $A_{\infty}$ structure and the construction of
Merkulov in \S 1.
In \S 2, we review the abstract Hodge theory of a differential graded algebra 
and show how it leads to the data in Merkulov's construction.
Applications to de Rham complex and deformation complex are presented
in \S 3.

\section{Merkulov's construction of construction of $A_{\infty}$ algebra}

Let $A = \oplus A^n$ be a $\bZ$-graded vector space.
An $A_{\infty}$ (algebra) structure on a vector space $A$
is a sequence of linear maps $m_k: A^{\otimes k} \to A$, $k \geq 1$,
$\deg m_k = 2 - k$
satisfying a sequence of conditions:

\begin{align*}
& m_1^2 = 0 \tag{A1}, \\
& m_1(m_2(a_1 \otimes a_2)) 
= m_2(m_1(a_1) \otimes a_2) + (-1)^{|a_1|} m_2(a_1 \otimes m_1(a_2)), 
	\tag{A2} \\
& m_2(m_2(a_1 \otimes a_2) \otimes a_3) 
- m_2(a_1 \otimes m_2(a_2 \otimes a_3)) \tag{A3} \\
&  = m_3(m_1(a_1) \otimes a_2 \otimes a_3) 
+ (-1)^? m_3(a_1 \otimes m_1(a_2) \otimes a_3) \notag \\
& + (-1)^? m_3(a_1 \otimes a_2 \otimes m_1(a_3)), \notag
\end{align*}
and so on. 
The formula ($A1$) says that $m_1$ is a differential.
If one writes $d = m_1$, then $(A, d)$ is a cochain complex.
Conversely, given a cochain complex, 
we can regard it as an $A_{\infty}$ algebra with $m _k = 0$ for $k \geq 2$.
Similarly, if we regard $m_2$ as a multiplication on $A$
and write $a_1 \cdot a_2 = m_2(a_1 \otimes a_2)$,
then ($A2$) says that $d$ is a derivation 
with respect to this multiplication.
Furthermore, if $m_3 = 0$,
then ($A3$) implies that the multiplication $\cdot$ is associative,
and hence $(A, \cdot, d)$ is a {\em differential graded algebra}.
Conversely, any differential graded algebra can be regarded as
an $A_{\infty}$ algebra with $m_k = 0$, $k \geq 3$.
In this paper, we will be interested in 
$A_{\infty}$ algebras with $m_1 = 0$,
which by ($A3$), 
are graded associative algebras with higher multiplications.

We now review Merkulov's construction.
Let $(V,d)$ be a differential graded algebra,
$W$ be a sub complex of $(V,d)$, i.e. vector subspace $W \subset V$
invariant under $d$.
$W$ is {\em not} assumed to be
a subalgebra of $V$. Instead assume that
there exists an odd operator
$Q: V \to V$
such that for any $v\in V$ the element $(1-[d,Q])v$ lies in the subspace
$W$, where $[\cdot , \cdot ]$ is the supercommutator.
Define a series of linear maps
$$
\lambda_n: \otimes^n V \to V, \ \ \ \ \ n\geq 2,
$$
starting with
$$
\lambda_2(v_1,v_2) := v_1\cdot v_2
$$
and then recursively, for $n\geq 3$,
\begin{eqnarray*}
&& \lambda_n(v_1,\ldots,v_n) \\
& = & (-1)^{n-1} [Q\lambda_{n-1}(v_1,\ldots,v_{n-1})]
\cdot v_n -
(-1)^{n|v_1|}v_1 \cdot [Q\lambda_{n-1}(v_2,\ldots, v_n)] \nonumber\\
&& -\,\sum_{\substack{k+l=n+1\\ k,l\geq 2}}(-1)^{k+(l-1)(|v_1|+\ldots+
|v_k|)} [Q\lambda_k(v_1, \ldots, v_k)] \cdot [Q\lambda_l(v_{k+1},\ldots,v_n)].
\nonumber \\
&& \label{la1}
\end{eqnarray*}

\begin{theorem} \label{main} (Merkulov \cite{Mer})
Let $(V,d)$ be a differential graded algebra and $(W,d)\subset
(V,d)$ be a subcomplex as above.
Then the linear maps
$$
m_k: \otimes^k W \to W, \ \ \  k\geq1,
$$
defined by
\begin{eqnarray}
m_1 &:=& d,\nonumber\\
m_k &:=& (1-[d,Q])\lambda_k, \ \ \mbox{for}\ {k\geq 2}, \label{mu}
\end{eqnarray}
with $\lambda_k$ being given above, satisfy the higher order
associativity identities.
I.e., $m_1, m_2, \cdots$ define an $A_{\infty}$ structure on $W$.
\end{theorem}

\section{Abstract Hodge theory}

Suppose that a differential graded algebra $(A, \wedge, \fd)$ is given
a (Euclidean or Hermitian) metric $\langle \cdot, \cdot \rangle$,
such that $\fd$ has a formal adjoint $\fd^*$,
i.e.,
$$\langle \fd a, b \rangle = \langle a, \fd^* b\rangle.$$
Since $\fd^2 = 0$,
it follows that $(\fd^*)^2 = 0$.
Set $\square_{\fd} = \fd\fd^* + \fd^* \fd$,
and $\cH = \Ker \square_{\fd} = \{a \in A: \square_{\fd} a = 0 \}$.
Then $\cH = \{a: \fd a = 0, \fd^* a = 0 \} = \Ker \fd \cap \Ker \fd^*$.
Assume that $A$ admits a ``Hodge decomposition'':
$$A = \cH \oplus \Img \fd \oplus \Img \fd^*.$$
It is standard to see that $\Ker \fd = \cH \oplus \Img \fd$,
and hence $H(A, \fd) \cong \cH$.
Notice that $\square_{\fd}|_{\Img \fd \oplus \Img \fd^*}$ is invertible.
Denote its inverse by $\square^{-1}$.
Consider the Green's operator $G: A \to A$, 
which is defined as the composition
$$A \to \Img \fd \oplus \Img \fd^* 
\stackrel{\square^{-1}}{\longrightarrow}
\Img \fd \oplus \Img \fd^*\to A,$$
where the first arrow is the inclusion,
and the last arrow is the projection.
Since $\square_{\fd}$ commutes with $\fd$ and $\fd^*$,
so does $G$.
For any $\alpha \in A$,
denote $f(\alpha)$ by $\alpha^H$. 
Then we have 
\begin{eqnarray} \label{decomposition}
\alpha - \alpha^H = \square_{\fd} G_{\fd} \alpha 
= (\fd\fd^* + \fd^*\fd) G_{\fd}\alpha
= \fd \fd^*G_{\fd}\alpha + \fd^* G_{\fd}\fd\alpha.
\end{eqnarray}
This gives the explicit Hodge decomposition.
Note that there is an induced wedge product on $H(A, \fd)$.

Hence $\cH$ has an induced product by virtue of being isomorphic to
$H(A, \fd)$ as vector spaces.
However, in general the product of two elements in $\cH$ may not lie in 
$\cH$.
For $\alpha, \beta \in \cH$, define
$$\alpha \circ \beta = (\alpha \wedge \beta)^H.$$
Merkulov \cite{Mer} made the remark that
this multiplication is in general not associative.
The following lemma shows that his remark is incorrect.

\begin{lemma}
For $\alpha, \beta, \gamma \in \cH$,
we have
$$((\alpha \wedge \beta)^H \wedge \gamma)^H 
= (\alpha \wedge \beta \wedge \gamma)^H.$$
Hence $(\cH, \circ)$ is associative.
\end{lemma}

\begin{proof}
By (\ref{decomposition}),
we have
\begin{eqnarray*}
 (\alpha \wedge \beta)^H 
& = & \alpha \wedge \beta - \fd \fd^*G_{\fd} (\alpha \wedge \beta)
+ \fd^* G_{\fd}\fd (\alpha \wedge \beta) \\
& = & \alpha \wedge \beta - \fd \fd^*G_{\fd} (\alpha \wedge \beta).
\end{eqnarray*}
Therefore,
\begin{eqnarray*}
(\alpha \wedge \beta)^H \wedge \gamma
& = & \alpha \wedge \beta \wedge \gamma 
 	- \fd \fd^*G_{\fd} (\alpha \wedge \beta) \wedge \gamma \\
& = & \alpha \wedge \beta \wedge \gamma 
 	- \fd (\fd^*G_{\fd} (\alpha \wedge \beta) \wedge \gamma).
\end{eqnarray*}
Hence
\begin{eqnarray*}
((\alpha \wedge \beta)^H \wedge \gamma)^H 
= (\alpha \wedge \beta \wedge \gamma)^H
- (\fd (\fd^*G_{\fd} (\alpha \wedge \beta) \wedge \gamma))^H
= (\alpha \wedge \beta \wedge \gamma)^H.
\end{eqnarray*}
As a consequence, 
\begin{eqnarray*}
 (\alpha \circ \beta) \circ \gamma 
& = & ((\alpha \wedge \beta)^H \wedge \gamma)^H 
= (\alpha \wedge \beta \wedge \gamma)^H \\
& = & (\alpha \wedge (\beta \wedge \gamma)^H)^H 
= \alpha \circ (\beta \circ \gamma).
\end{eqnarray*}
\end{proof}

In fact, another to check the associativity is the following

\begin{lemma} \label{isomorphism}
The isomorphism $\phi: \cH \to H(A, \fd)$ given by $\alpha \mapsto [\alpha]$
(where $[\alpha]$ denote the cohomology class of $\alpha$)
is a ring isomorphism $(\cH, \circ) \to (H(A, \fd), \wedge)$.
\end{lemma}

\begin{proof}
For $\alpha, \beta \in \cH$,
we have
\begin{eqnarray*}
\phi( \alpha \circ \beta) = [(\alpha \wedge \beta)^H] 
= [\alpha \wedge \beta] = [\alpha] \wedge [\beta] 
= \phi(\alpha) \wedge \phi(\beta).
\end{eqnarray*}
\end{proof}

Our main observation is that 
it is possible construct an $A_{\infty}$ structure on $\cH$
by the method of Merkulov \cite{Mer} reviewed in last section.
Indeed, let $V = A$, $W = \cH$, 
and $Q = G \fd^*$. 
Then (\ref{decomposition}) implies that
$(W \subset V, \fd, Q)$ satisfies the conditions in 
last section.
By Theorem \ref{main},
we get

\begin{theorem} \label{thm:Ainfty1}
For a differential graded algebra $(A, \wedge, \fd)$ 
with a Euclidean or Hermitian metric,
such that $\fd$ has a formal adjoint $\fd^*$ and $A$ has a Hodge decomposition
$$A = \cH \oplus \Img \fd \oplus \Img \fd^*,$$
where $\cH = \Ker \square_{\fd}$, 
there is a canonical $A_{\infty}$ structure on $\cH$
with $m_2 = \circ$.
\end{theorem}

By Lemma \ref{isomorphism},
we get 

\begin{theorem} \label{thm:Ainfty2}
For a differential graded algebra $(A, \wedge, \fd)$ 
with a Euclidean or Hermitian metric,
such that $\fd$ has a formal adjoint $\fd^*$ and $A$ has a Hodge decomposition
$$A = \cH \oplus \Img \fd \oplus \Img \fd^*,$$
there is a canonical $A_{\infty}$ structure on $H(A, \fd)$
with $m_2$ the induced wedge product.
\end{theorem}

\section{Applications}

\subsection{The Riemannian case}
Let $(M, g)$ be an oriented closed Riemannian manifold.
Then there is an induced Euclidean metric on $\Omega^*(M)$.
Let $\Omega^*(M)$ be the space of differential forms on $M$,
it has a wedge product $\wedge$ and the exterior differential $d$,
such that  $(\Omega^*(M), \wedge, d)$ is
a differential graded algebra.
We call it the {\em de Rham algebra}.
Using the orientation, one can define the Hodge star operator $*$.
The formal adjoint of $d$ is given by
$$d^* = - * d *.$$
Let $\square = dd^* + d^* d$.
Then $\square$ is a second order elliptic operator,
and standard elliptic operator theory gives the Hodge decomposition
$$\Omega^*(X) = \cH \oplus \Img d \oplus \Img d^*,$$
where $\cH$ is the space of harmonic forms.
Theorem \ref{thm:Ainfty2} then yields a canonical $A_{\infty}$ structure on 
the de Rham cohomology $H^*(M) = H^*(\Omega^*(M), d)$.

\subsection{The oomplex case}

Let $W$ be a closed complex manifold,
\begin{eqnarray*}
&& \Omega^{*, *}(W) = \oplus_{p, q \geq 0} 
\Gamma(W, \Lambda^p T^*W \otimes \Lambda^q \overline{T}^*W), \\
&& \Omega^{-*, *}(W) = \oplus_{p, q \geq 0} 
\Gamma(W, \Lambda^p TW \otimes \Lambda^q \overline{T}^*W).
\end{eqnarray*}
On both of these spaces, 
there is a wedge product $\wedge$ and a $\bar{\partial}$ operator,
which form differential graded algebras.
We call them {\em Dolbeault algebra} and {\em deformation algebra}
respectively.
Their cohomologies are denoted by 
$H^{*, *}(W)$ and $H^{-*, *}(W)$ respectively.
Given a Hermitian metric on the holomorphic tangent bundle $TW$,
there are induced Hermitian metrics on both
$\Omega^{*, *}(W)$ and $\Omega^{-*, *}(W)$.
It is easy to find the formal adjoints of the operators $\bar{\partial}$
on these spaces and define the corresponding Laplacian operators.
Again we get elliptic operators and by standard theory 
the Hodge decompositions.
Therefore, we have canonical $A_{\infty}$ structures on both 
$H^{*, *}(W)$ and $H^{-*, *}(W)$.

\subsection{The bundle case}
Let $X$ be a closed complex manifold,
and $\pi: E \to X$ a holomorphic vector bundle.
Consider the holomorphic vector bundle $\End(E) \to X$
and the space $\cA = \Omega^{0, *}(\End(E)
= \Gamma(X, \Lambda^*\overline{T}^*X \otimes \End(E))$.
There is an induced wedge product on $\cA$,
which is in general not graded commutative.
The $\bar{\partial}_{\End(E)}$ operator for $\End(E)$ gives 
a differential for $(\cA, \wedge)$.
Given a Hermitian metric on $X$ and a Hermitian metric on $E$,
one can consider 
the formal adjoint operator $\bar{\partial}_{\End(E)}^*$.
The corresponding Laplacian operator is elliptic,
hence we have Hodge theory.
Therefore, we have an induced $A_{\infty}$ structure on the Dolbeault
cohomology $H^*(X, \End(E))$.

\medskip

{\bf Acknowledgement}.
{\em The work in paper is done while 
the author is visiting Texas A$\&$M University.
It is partially supported by an NSF group infrastructure grant.
The author thanks the Mathematics Department and the Geometry-Analysis-Topology
group for hospitality and financial support.
Special thanks are due to Huai-Dong Cao. 
The collaboration with him expands the author's scope of interests and
brings his attention to the problem considered in this paper.
Interactions with Jim Stasheff enhance the author's interest and understanding
of infinite structures and many other things.}

\end{document}